\documentclass[a4paper]{article}
\usepackage{amssymb}

\newtheorem{theorem}{Theorem}

\newtheorem{corollary}[theorem]{Corollary}

\newtheorem{lemma}[theorem]{Lemma}

\newtheorem{proposition}[theorem]{Proposition}

\input{tcilatex}
\begin{document}

\begin{center}
\textbf{Several Metric Properties of Level Curves}

\medskip

Pisheng Ding

\medskip

\textsc{January, 2017}

\medskip
\end{center}

\begin{quote}
\textsc{Abstract. }{\small This article establishes several remarkably
simple identities relating certain metric invariants of level curves of real
and complex functions. In particular, we relate lengths of level curves to
their curvature and to the gradient field of the function. Some geometric
and analytic applications of the results are shown.}
\end{quote}

\section{Introduction}

This article examines certain metric invariants of level curves of a real
two-variable function and of a complex analytic function. In particular, we
relate lengths of level curves to their curvature and to the gradient field.

We review some terminologies before stating the main result. A \textit{Morse
function} on an open subset of $%
\mathbb{R}
^{n}$ is a real-valued $C^{2}$ function whose critical points are
nondegenerate (i.e., with nonsingular Hessian). Near a critical point, a
Morse function behaves as a nondegenerate quadratic form (i.e., its
second-degree Taylor approximation) and thus has no other critical points. A 
\textit{regular value} of a Morse function is a number which is not the
image of a critical point\footnote{%
As a seemingly peculiar (but useful for certain differential-topological
purposes) consequence of the definition, any number \textit{not} in the
range of the function \textit{is} a regular value. However, we will only be
concerned with regular values that are attained by the function.}, whereas a 
\textit{critical value} is the image of a critical point.

Let $f$ be a Morse function on an open connected subset of $%
\mathbb{R}
^{2}$. Suppose that $a$ and $b$ are values of $f$ such\ that $f^{-1}([a,b])$
is compact. For $t\in \lbrack a,b]$, let $L(t)$ be the length of the level-$%
t $ curve $f^{-1}(t)$; $L(t)$ is well-defined even if $t$ is a critical
value. At each \textit{regular point} (i.e., noncritical point) on $%
f^{-1}(t) $, install a unit tangent $\mathbf{T}$ and a unit normal $\mathbf{N%
}$ by letting $\mathbf{N}=-\nabla f/\left\vert \nabla f\right\vert $ and
requiring that the frame $(\mathbf{T},\mathbf{N})$ be positively-oriented.
Then, the \textit{signed}\ curvature $\kappa $ of $f^{-1}(t)$ at all regular
points is defined by the equation $d\mathbf{T}/ds=\kappa \mathbf{N}$, where $%
s$ is arc length (with its positive direction induced by $\mathbf{T}$).
Thereby, $\kappa $ is defined at all regular points of $f^{-1}([a,b])$.

With the notations introduced, we now state our main result.

\smallskip

\noindent \textbf{Main Result}:\quad \textit{If }$a$\textit{\ and }$b$%
\textit{\ are values of }$f$\textit{\ such\ that }$f^{-1}([a,b])$\textit{\
is compact, then, for }$t\in \lbrack a,b]$\textit{,}%
\[
\frac{d}{dt}\iint\nolimits_{f^{-1}([a,t])}\left\vert \nabla f\right\vert
dA\;=\;L(t)\;=\;L(a)+\iint\nolimits_{f^{-1}([a,t])}\kappa \,dA\,\text{.} 
\]

\noindent By $\iint\nolimits_{f^{-1}([a,t])}\kappa \,dA$, we mean either the
(proper) Riemann integral of $\kappa $ when $f^{-1}([a,t])$ is free of
critical points, or the improper Riemann integral of $\kappa $ otherwise;
the latter case will be addressed in detail.

These two identities are established in \S 2 and their counterparts for
complex analytic functions are given in \S 3. A number of applications are
shown in \S 4, including in \S 4.1 instances of curve evolution and a
characterization of circles and in \S 4.2 parallel results for level
surfaces. Two technical details are relegated to \S 5.

\section{Level Curves of Real Functions}

\subsection{Preliminaries}

Let $\Omega $ be an open connected subset of $%
\mathbb{R}
^{2}$ and $f$ be a Morse function (defined in \S 1) on $\Omega $. The
notions of \textit{regular value}, \textit{critical value}, and \textit{%
regular point} of $f$ have been introduced in \S 1.

We shall mean, by a \textit{regular }(resp. \textit{critical})\textit{\ level%
} of $f$, a nonempty level set $f^{-1}(t)$ with $t$ a regular (resp.
critical) value. Near a regular point $P$, the set $f^{-1}(f(P))$ is locally
a $C^{2}$ curve (by the implicit function theorem). Therefore, each
component of a regular level is globally a $C^{2}$ curve, and a compact
regular level is a disjoint union of simple closed curves.

Let $[a,b]$ be an interval of regular values attained by $f$ such that $%
f^{-1}([a,b])$ is compact. As basic facts in differential topology, $%
f^{-1}(t)$ is diffeomorphic to $f^{-1}(a)$ for $t\in \lbrack a,b]$ and $%
f^{-1}([a,b])$ is diffeomorphic to $f^{-1}(a)\times \lbrack a,b]$. The idea
is to start a flow originating from $f^{-1}(a)$ following the gradient field
and to control the flow's speed so that points on $f^{-1}(a)$ all
\textquotedblleft drift\textquotedblright\ to points of the same $f$-value
at the same time; details are given in \S 2.2.

Suppose that $f^{-1}(t)$ is compact and consider its length $L(t)$. We
explain that $L(t)$ is well-defined even if $t$ is a critical value. Near a
critical point $P\in f^{-1}(t)$, a local $C^{2}$ coordinate system $%
(u(x,y),v(x,y))$ exists such that $u(P)=0=v(P)$ and $f(x,y)=t+q(u,v)$ where $%
q$ is a nondegenerate quadratic form; see \cite[pp.\thinspace 54-58]{Wallace}%
. A neighborhood of $P$ on $f^{-1}(t)$ is thus diffeomorphic to a
neighborhood of the origin on $q^{-1}(0)$, which is sufficiently
well-behaved to admit length. If $f^{-1}([a,b])$ is compact, then $L$ is
continuous on $[a,b]$; its continuity at a regular value $t_{0}$ is due to
the diffeomorphism between $f^{-1}([t_{0}-\epsilon ,t_{0}+\epsilon ])$ and $%
f^{-1}(t_{0})\times \lbrack t_{0}-\epsilon ,t_{0}+\epsilon ]$, whereas its
continuity at a critical value is due to the good behavior of $q^{-1}(\delta
)$ as $\delta $ varies over $(-\epsilon ,\epsilon )$.

Recall from \S 1 our choice of $\mathbf{N}$ and $\mathbf{T}$ at regular
points of $f^{-1}(t)$; i.e., $\mathbf{N}:=-\nabla f/\left\vert \nabla
f\right\vert $ and $(\mathbf{T},\mathbf{N})$ is a positively-oriented frame.
At regular points on $f^{-1}(t)$, the \textit{signed} curvature $\kappa $
and the frame $(\mathbf{T},\mathbf{N})$ are (by definition) related by%
\[
\frac{d\mathbf{T}}{ds}=\kappa \mathbf{N}\text{\quad and\quad }\frac{d\mathbf{%
N}}{ds}=-\kappa \mathbf{T\,}\text{,} 
\]%
where $s$ is arc length (with its positive direction induced by $\mathbf{T}$%
). For a simple closed component $C$ of a regular level,%
\[
\int_{C}\kappa \,ds=\pm 2\pi 
\]%
where the plus sign is in force iff $\mathbf{T}$ induces the positive
orientation\footnote{%
The positive orientation on a simple closed curve $C$ is the one following
which the traversal of $C$ gives a positive winding number around any
interior point in the Jordan domain bounded by $C$; see \cite[pp.\thinspace
392-396]{do Carmo}. If $C$ is convex, the positive orientation on $C$ is
counterclockwise.} on $C$; see \cite[pp.\thinspace 36-37]{do Carmo}.

Henceforth, we will regard $\kappa $ as a function on the set of regular
points in $\Omega $, i.e., $\kappa (P)$ is the signed curvature at $P$ of
the curve $f^{-1}(f(P))$.

It is a fact (shown in \cite[p.\thinspace 125]{Courant}) that%
\begin{equation}
\kappa =\frac{f_{xx}f_{y}^{2}-2f_{xy}f_{x}f_{y}+f_{yy}f_{x}^{2}}{\left\vert
\nabla f\right\vert ^{3}}\text{ ;}  \label{EQ Curvature}
\end{equation}%
it is a matter of computation to verify that%
\begin{equation}
\func{div}\left( \frac{\nabla f}{\left\vert \nabla f\right\vert }\right)
=\kappa \,\text{.}  \label{EQ Divergence}
\end{equation}%
As $\kappa $ becomes unbounded near a critical point $P$, its improper
Riemann integral $\iint\nolimits_{D}\kappa \,dA$ over a small disc $D$
around $P$ nonetheless converges. We relegate this matter to \S 5, in which
we prove the integrability of $\kappa $ under a hypothesis that generalizes
the Morse condition; before then, we take this fact for granted.

\subsection{Length in Relation to Curvature and Gradient}

For a number of results, we assume the following hypothesis.

\smallskip

\noindent \textbf{Hypothesis (\dag )}$\quad f$ is a $C^{2}$ Morse function
on an open connected subset $\Omega $ of $%
\mathbb{R}
^{2}$; $a$ and $b$ are values attained by $f$ such that $f^{-1}([a,b])$ is
compact. (As $\Omega $ is connected, each $t\in \lbrack a,b]$ is also a
value attained by $f$.)

\smallskip

We begin with a key lemma, upon which many of our results rest.

\begin{lemma}
\label{Lemma KEY}Assume Hypothesis (\dag ). If a function $g$ is continuous
on the set of regular points in $f^{-1}([a,b])$ and is Riemann-integrable on 
$f^{-1}([a,b])$, then%
\begin{equation}
\iint_{f^{-1}([a,b])}g\,dA=\int_{a}^{b}\left( \int_{f^{-1}(t)}\frac{g}{%
\left\vert \nabla f\right\vert }ds\right) dt\,\text{,}  \nonumber
\end{equation}%
where the line integral $\int_{f^{-1}(t)}\left( g/\left\vert \nabla
f\right\vert \right) ds$ is only defined for $t$ a regular value.\footnote{%
The improper Riemann integral $\int_{a}^{b}\varphi (t)dt$ is a well-defined
notion when $\varphi $ is continuous at all but finitely many $t\in \lbrack
a,b]$. In the present case, $\varphi (t)=\int_{f^{-1}(t)}\left( g/\left\vert
\nabla f\right\vert \right) ds$.}
\end{lemma}

When $[a,b]$ is free of critical values, this formula is (somewhat
informally) shown in \cite[p.\thinspace 298-300]{Courant}. We shall sketch a
proof based on a differential-topological construction, which will be of use
in a later argument.

\medskip

\noindent \textit{Proof. }There are two cases, according as whether $[a,b]$
contains a critical value.

\noindent \textit{Case 1.}\quad Suppose that $[a,b]$ contains no critical
value. For each $p\in f^{-1}(a)$, let $t\mapsto H(p,t)$ be the flow (i.e.
the integral curve) for the field $\nabla f/\left\vert \nabla f\right\vert
^{2}$, originating from $p$ at the initial time $t=a$; i.e., $H(p,-)$ is the
solution to the initial value problem%
\begin{equation}
\frac{\partial H}{\partial t}(p,t)=\frac{\nabla f\left( H(p,t)\right) }{%
\left\vert \nabla f\left( H(p,t)\right) \right\vert ^{2}}\text{\quad
with\quad }H(p,a)=p\,\text{.}  \label{EQ H}
\end{equation}%
Clearly, $f(H(p,t))=t$, as $f(H(p,t))-f(H(p,a))=\int_{a}^{t}\nabla f\cdot
\left( \partial H/\partial \tau \right) d\tau =\int_{a}^{t}1d\tau $. Due to
existence, uniqueness, and smooth dependence on initial condition of
solution to the initial value problem, the map $H(-,t):f^{-1}(a)\rightarrow
f^{-1}(t)$, and consequently the map $H:f^{-1}(a)\times \lbrack
a,b]\rightarrow f^{-1}([a,b])$, are diffeomorphisms; see \cite[p.\thinspace
65]{Wallace} and the Remark that follows this proof.

Now assume that $f^{-1}(a)$ is connected and hence a simple closed curve;
otherwise, treat each component of $f^{-1}(a)$ and its corresponding
component of $f^{-1}([a,b])$ separately. Let $\gamma (-,a):[0,1]\rightarrow 
\mathbb{R}
^{2}$ be a regular parametrization of $f^{-1}(a)$ that is one-to-one on $%
[0,1)$ with $\gamma (0,a)=\gamma (1,a)$ and that induces the same
orientation that $\mathbf{T}$ does on $f^{-1}(a)$. The map%
\begin{equation}
\gamma (u,t):=H(\gamma (u,a),t)  \label{EQ Parametrization}
\end{equation}%
then parametrizes $f^{-1}([a,b])$, with%
\[
\left\vert \det \frac{\partial (x,y)}{\partial (u,t)}\right\vert =\left\vert 
\frac{\partial \gamma }{\partial t}\times \frac{\partial \gamma }{\partial u}%
\right\vert =\left\vert \frac{\partial \gamma }{\partial t}\right\vert
\left\vert \frac{\partial \gamma }{\partial u}\right\vert =\frac{1}{%
\left\vert \nabla f\right\vert }\left\vert \frac{\partial \gamma }{\partial u%
}\right\vert \,\text{.} 
\]%
Hence, by changing variables in integration,%
\[
\iint_{f^{-1}([a,b])}g\,dA=\int_{a}^{b}\left( \int_{0}^{1}\frac{g}{%
\left\vert \nabla f\right\vert }\cdot \left\vert \frac{\partial \gamma }{%
\partial u}\right\vert du\right) dt=\int_{a}^{b}\left( \int_{f^{-1}(t)}\frac{%
g}{\left\vert \nabla f\right\vert }ds\right) dt\,\text{.} 
\]%
\textit{Case 2.}\quad Suppose that $[a,b]$ contains critical values
comprising the (finite) set $S$. Then, $(a,b)\smallsetminus S$ is a disjoint
union of finitely many intervals $I_{j}:=(c_{j},c_{j+1})$ of regular values
attained by $f$. As $f^{-1}([a,b])=\cup _{j}f^{-1}(I_{j})\cup f^{-1}\left(
S\cup \{a,b\}\right) $ and $f^{-1}\left( S\cup \{a,b\}\right) $ has zero
area,%
\[
\iint_{f^{-1}([a,b])}g\,dA=\sum_{j}\iint_{f^{-1}(I_{j})}g\,dA\,\text{.} 
\]%
Applying Case 1 to $f^{-1}([c_{j}+\epsilon ,c_{j+1}-\delta ])$ and letting $%
\epsilon ,\delta \rightarrow 0^{+}$, we have%
\begin{eqnarray*}
\iint_{f^{-1}(I_{j})}g\,dA &=&\lim_{\epsilon ,\delta \rightarrow
0}\iint_{f^{-1}([c_{j}+\epsilon ,c_{j+1}-\delta ])}g\,dA \\
&=&\lim_{\epsilon ,\delta \rightarrow 0}\int_{c_{j}+\epsilon
}^{c_{j+1}-\delta }\left( \int_{f^{-1}(t)}\frac{g}{\left\vert \nabla
f\right\vert }ds\right) dt \\
&=&\int_{I_{j}}\left( \int_{f^{-1}(t)}\frac{g}{\left\vert \nabla
f\right\vert }ds\right) dt\,\text{.}
\end{eqnarray*}%
Summing these integrals over $j$ proves the assertion.\hfill $\blacksquare $

\medskip

Choices for $g$ in Lemma \ref{Lemma KEY} yield integral identities. For
example, the choice $g\equiv 1$ gives the area of $f^{-1}([a,b])$ as $%
\int_{a}^{b}\left( \int_{f^{-1}(t)}1/\left\vert \nabla f\right\vert
ds\right) dt$. We next show that Lemma \ref{Lemma KEY} can result in
interesting relations among the length of level curves, their curvature, and
the gradient field.

\begin{theorem}
\label{Theorem MAIN}Assume Hypothesis (\dag ).

\begin{description}
\item[(a)] $\iint_{f^{-1}[a,b]}(h\circ f)\cdot \left\vert \nabla
f\right\vert dA=\int_{a}^{b}h(t)L(t)\,dt$ for any function $h$ such that
both members of the equation are meaningful. In particular, for any $t\in
\lbrack a,b]$,%
\[
\int_{a}^{t}L(\tau )\,d\tau =\iint_{f^{-1}([a,t])}\left\vert \nabla
f\right\vert dA\,\text{,} 
\]%
or equivalently, 
\[
L(t)=\frac{d}{dt}\iint_{f^{-1}([a,t])}\left\vert \nabla f\right\vert dA\,%
\text{.} 
\]

\item[(b)] For any regular value $t\in \lbrack a,b]$,%
\[
L^{\prime }(t)=\int_{f^{-1}(t)}\frac{\kappa }{\left\vert \nabla f\right\vert 
}ds\,\text{,} 
\]%
or equivalently, for any $t\in \lbrack a,b]$,%
\[
L(t)=L(a)+\iint\nolimits_{f^{-1}([a,t])}\kappa \,dA\,\text{.} 
\]
\end{description}
\end{theorem}

As promised in \S 2.1, the integrability of $\kappa $ will be proven in \S 5.

\medskip

\noindent \textit{Proof. }Part (a) follows from Lemma \ref{Lemma KEY} by
letting $g=(h\circ f)\cdot \left\vert \nabla f\right\vert $; the
\textquotedblleft particular\textquotedblright\ case results from letting $%
h=\chi _{\left[ a,t\right] }$, the characteristic function of the interval $%
[a,t]$. (Continuity of $L$ makes applicable the fundamental theorem of
calculus, i.e., $\frac{d}{dt}\int_{a}^{t}L(\tau )d\tau =L(t)$.)

In Part (b), the equivalence between the two assertions is due to Lemma \ref%
{Lemma KEY}, as we now show. First, assume the formula for $L^{\prime }$.
Since $L$ is continuous at the (finitely many) critical values in $[a,b]$,
the fundamental theorem of calculus applies to give, for $t\in \lbrack a,b]$,%
\[
L(t)-L(a)=\int_{a}^{t}L^{\prime }(\tau )\,d\tau =\int_{a}^{t}\left(
\int_{f^{-1}(\tau )}\frac{\kappa }{\left\vert \nabla f\right\vert }ds\right)
d\tau =\iint\nolimits_{f^{-1}([a,t])}\kappa \,dA\,\text{,} 
\]%
where the last equality follows from Lemma \ref{Lemma KEY} (with $t$ playing
the role of $b$). Now assume the integral formula for $L$. Then, for a
regular value $t\in (a,b)$,%
\[
L^{\prime }(t)=\frac{d}{dt}\iint\nolimits_{f^{-1}([a,t])}\kappa \,dA=\frac{d%
}{dt}\int_{a}^{t}\left( \int_{f^{-1}(\tau )}\frac{\kappa }{\left\vert \nabla
f\right\vert }ds\right) d\tau =\int_{f^{-1}(t)}\frac{\kappa }{\left\vert
\nabla f\right\vert }ds\text{\thinspace ,} 
\]%
where, in the last equality, $t$ being a regular value is essential, as $%
\frac{d}{dt}\int_{a}^{t}\varphi (\tau )d\tau =\varphi (t)$ iff $\varphi $ is
continuous at $t$.

We now give two alternative proofs for Part (b) by directly proving each of
the two equivalent formulae.

To prove the formula for $L^{\prime }$, let $t\in (a,b)$ be a regular value.
There is an interval $[a^{\prime },b^{\prime }]$ of regular values such that 
$t\in (a^{\prime },b^{\prime })\subset \lbrack a,b]$. Let $a^{\prime }$ and $%
b^{\prime }$ play the role of $a$ and $b$ in case 1 of the proof of Lemma %
\ref{Lemma KEY} and let%
\[
\gamma :[0,1]\times \lbrack a^{\prime },b^{\prime }]\rightarrow
f^{-1}([a^{\prime },b^{\prime }]) 
\]%
be as in (\ref{EQ Parametrization}). For the speed $v:=\left\vert \partial
\gamma /\partial u\right\vert $, we have $v^{2}=\left\langle \partial \gamma
/\partial u,\partial \gamma /\partial u\right\rangle $ and%
\begin{eqnarray*}
\frac{1}{2}\frac{\partial (v^{2})}{\partial t} &=&\left\langle \frac{%
\partial }{\partial u}\left( \frac{\partial \gamma }{\partial t}\right) ,%
\frac{\partial \gamma }{\partial u}\right\rangle =\left\langle \frac{%
\partial }{\partial u}\left( \frac{\nabla f}{\left\vert \nabla f\right\vert
^{2}}\right) ,\frac{\partial \gamma }{\partial u}\right\rangle \\
&=&\left\langle \frac{\partial }{\partial u}\left( \frac{-1}{\left\vert
\nabla f\right\vert }\mathbf{N}\right) ,v\mathbf{T}\right\rangle
=\left\langle \frac{-1}{\left\vert \nabla f\right\vert }\frac{\partial 
\mathbf{N}}{\partial u},v\mathbf{T}\right\rangle \text{\thinspace .}
\end{eqnarray*}%
As $\partial \mathbf{N}/\partial u=v\partial \mathbf{N}/\partial s=-v\kappa 
\mathbf{T}$,%
\[
\frac{\partial v}{\partial t}=\frac{1}{2v}\frac{\partial (v^{2})}{\partial t}%
=\frac{1}{v}\left\langle \frac{1}{\left\vert \nabla f\right\vert }v\kappa 
\mathbf{T},v\mathbf{T}\right\rangle =\frac{\kappa v}{\left\vert \nabla
f\right\vert }\text{\thinspace .} 
\]%
and%
\[
L^{\prime }(t)=\frac{d}{dt}\int_{0}^{1}v\,du=\int_{0}^{1}\frac{\partial v}{%
\partial t}du=\int_{0}^{1}\frac{\kappa }{\left\vert \nabla f\right\vert }%
vdu=\int_{f^{-1}(t)}\frac{\kappa }{\left\vert \nabla f\right\vert }ds\,\text{%
.} 
\]%
\hfill

To prove the integral formula for $L$, first consider the case when $[a,b]$
is free of critical values. Let $R$ denote $f^{-1}([a,b])$ in this argument.
Let $\mathbf{n}$ denote the unit outward normal (relative to $R$) on $%
\partial R$; it should be clear that $\mathbf{n}=-\nabla f/\left\vert \nabla
f\right\vert $ on $f^{-1}(a)$ and $\mathbf{n}=\nabla f/\left\vert \nabla
f\right\vert $ on $f^{-1}(b)$. Then, by Green's theorem and (\ref{EQ
Divergence}),%
\[
L(b)-L(a)=\int_{\partial R}\left\langle \frac{\nabla f}{\left\vert \nabla
f\right\vert },\mathbf{n}\right\rangle ds=\iint\nolimits_{R}\func{div}\left( 
\frac{\nabla f}{\left\vert \nabla f\right\vert }\right)
dA=\iint\nolimits_{R}\kappa \,dA\,\text{.} 
\]%
If $[a,b]$ has critical values only in its interior, this argument can be
modified by excising from $R$ small discs around the critical points (where $%
\nabla f/\left\vert \nabla f\right\vert $ ceases to be defined); if $a$
(resp. $b$) is a critical value, replace it by $a+\epsilon $ (resp. $%
b-\epsilon $). We omit the details but note that the integrability of $%
\kappa $ over $R$ is essential.$\hfill \blacksquare $

\medskip

Comparing the two proofs given for Theorem \ref{Theorem MAIN}(b), the one
deriving the formula for $L^{\prime }$ is perhaps more geometrically
revealing, for, besides its independence on (\ref{EQ Divergence}), it shows
directly how curvature and gradient affect length of level curves. Indeed,
varying metric properties of a family of curves belong to the broader field
of evolution of curves and surfaces, which encompasses topics such as the
curve-shortening flow and mean curvature flow. In \S 4.1, we show how our
results on level curves may inform certain instances of curve evolution.

\medskip

\noindent \textbf{Remark}\textit{\quad }Assume Hypothesis (\dag ). By
Theorem \ref{Theorem MAIN}(a), we have the estimate%
\[
\min_{t\in \lbrack a,b]}L(t)<\frac{1}{b-a}\iint_{f^{-1}([a,b])}\left\vert
\nabla f\right\vert dA<\max_{t\in \lbrack a,b]}L(t)\,\text{.} 
\]%
This inequality is optimal over all $f$ satisfying the hypothesis, as
evidenced by the functions $r^{n}$ ($n\in 
\mathbb{N}
$) on the unit disc.$\hfill \blacksquare $

\medskip

\noindent \textbf{Remark}\quad Suppose that $R\subset \Omega $ is a closed
Jordan domain bounded by a component of a level set of $f$. Then, the length 
$\mathcal{L}(\partial R)$ of $\partial R$ has a simple formula%
\begin{equation}
\mathcal{L}(\partial R)=\pm \iint\nolimits_{R}\kappa \,dA\,\text{.}
\label{EQ Jordan Boundary}
\end{equation}%
The proof for the integral formula in Theorem \ref{Theorem MAIN}(b) will
also prove this result, once we note that no critical points of $f$ lie on $%
\partial R$ (due to our assumption that $f$ is a Morse function\footnote{%
While a level curve passing through a saddle point of a \textit{Morse
function} necessarily self-intersects, this may not be so in general, as
evidenced by the level-$0$ curve of the function $g(x,y)=y^{3}-x^{6}$.});
the integrability of $\kappa $ is again essential, as $R$ necessarily
contains critical points in its interior.$\hfill \blacksquare $

\section{Level Curves of Complex Analytic Functions}

\subsection{Level Curves near a Critical Point}

We explain why level curves of a complex analytic function are \textit{%
locally} so well-behaved that their intersections with a small disc around
any point admit continuous length.

For $z_{0}\in 
\mathbb{C}
$ and $r>0$, let%
\[
D(z_{0};r)=\left\{ z:\left\vert z-z_{0}\right\vert <r\right\} \text{\quad
and\quad }C(z_{0};r)=\left\{ z:\left\vert z-z_{0}\right\vert =r\right\} \,%
\text{.} 
\]%
Let $f$ denote a nonconstant analytic function on an open connected subset $%
\Omega $ of $%
\mathbb{C}
$. By the\textit{\ level-}$t$\textit{\ set of }$f$, we mean the level-$t$
set of $\left\vert f\right\vert $. As $\left\vert f\right\vert $ is not
differentiable at any simple zero of $f$, we will consider $\left\vert
f\right\vert ^{2}$ when classifying level curves. By a \textit{regular
(resp. critical)\ level} of $f$, we mean a regular (resp. critical) level of 
$\left\vert f\right\vert ^{2}$.\ Let $L(t)$ denote the length (whenever it
can be defined)\ of the level-$t$ set of $f$.

A simple calculation using the Cauchy-Riemann equations shows that%
\[
\left\vert \nabla \left( \left\vert f\right\vert ^{2}\right) \right\vert
=2\left\vert f\cdot f^{\prime }\right\vert \,\text{.} 
\]%
The critical points of $\left\vert f\right\vert ^{2}$ are then the zeros of $%
f$ and of $f^{\prime }$. If $|f(z_{0})|^{2}=t^{2}$ ($t>0$) is a regular
value of $\left\vert f\right\vert ^{2}$, then, since $f^{\prime }(z_{0})\neq
0$, $f$ conformally maps an arc $C\ni z_{0}$ on the level-$t$ curve onto an
arc on the circle $C(0;t)$. In other words, $C$ is the conformal image
(under the inverse of $f|_{D(z_{0};\delta )}$) of a circular arc; hence a
regular level of $f$ is an analytic curve, the smoothest of all.

However, $\left\vert f\right\vert ^{2}$ may not be a Morse function; its
level curves near a saddle point are somewhat more complicated but still
have a simple description. For simplicity, we will call $z_{0}$ a \textit{%
saddle point} of $f$ if it is a saddle point of $\left\vert f\right\vert
^{2} $. If $z_{0}$ a saddle point of $f$, then $f(z_{0})\neq 0$ but $%
f^{\prime }(z_{0})=0$, and we show that level curves of $f$ within $%
D(z_{0};\epsilon )$ are images under a conformal map\ of level curves of $%
z^{n}-1$ within $D(0;\delta )$, where $n$ is the order of $z_{0}$ as a zero
of $f(z)-f(z_{0})$. (For contrast, recall that the behavior of a nondefinite
quadratic form near the origin dictates that of a Morse function near a
saddle point.)

In general, if $f$ is analytic at $z_{0}$, then there is a neighborhood $U$
of $z_{0}$ and a conformal equivalence $w:U\rightarrow D(0;\delta )$ with $%
\delta <1$ such that%
\[
f(z)=\left[ w(z)\right] ^{n}+f(z_{0})\quad \text{for\quad }z\in U\,\text{;} 
\]%
necessarily, $w(z_{0})=0$ and $n=\min \left\{ k:f^{(k)}(z_{0})\neq 0\right\} 
$.

Suppose now that $z_{0}$ is a critical point of $\left\vert f\right\vert
^{2} $, in which case either $f(z_{0})=0$ or $f^{\prime }(z_{0})=0$.

If $f(z_{0})=0$, then the level-$0$ set of $f|_{U}$ is the singleton $%
\{z_{0}\}$, while the level-$\epsilon $ set of $f|_{U}$ (for small $\epsilon
>0$) is the level-$\sqrt[n]{\epsilon }$ set of $w$, which is the (conformal)
image under $w^{-1}$ of the origin-centered circle of radius $\sqrt[n]{%
\epsilon }$. So the length $L(t)$ of the level-$t$ curves of $f|_{U}$ is
continuous at $0=\left\vert f(z_{0})\right\vert $.

If $f(z_{0})\neq 0$ but $f^{\prime }(z_{0})=0$, then, by considering $%
-f(z)/f(z_{0})$, we may assume that $f(z_{0})=-1$, in which case $f=w^{n}-1$
on $U$ (with $n\geq 2$). For $t$ near the critical value $1$, the level-$t$
set of $f|_{U}$ is the image under $w^{-1}$ of the set $\left\{ \zeta \in
D(0;\delta ):\left\vert \zeta ^{n}-1\right\vert =t\right\} $. Hence, it
suffices to consider the level-$t$ curves (with $t$ near $1$) of the
function $p_{n}:\zeta \mapsto \left( \zeta ^{n}-1\right) $ on $D(0;\delta )$%
, to which we now turn.

Denote the input for $p_{n}$ by $z$ (instead of $\zeta $). With $\Gamma _{t}$
denoting the level-$t$ set of $p_{n}$ in $D(0;\delta )$, i.e.,%
\[
\Gamma _{t}:=\left\{ z\in D(0;\delta ):\left\vert z^{n}-1\right\vert
=t\right\} \text{,}\, 
\]%
we have%
\[
z\in \Gamma _{t}\text{ \ iff \ }z^{n}\in \Sigma _{t}:=C(1;t)\cap D(0;\delta
^{n})\text{ ;} 
\]%
in other words, $\Gamma _{t}$ is comprised of all the $n$th roots of every
complex number on the circular arc $\Sigma _{t}$.

As $t$ ranges over $[1-\epsilon ,1]$, the circular arc $\Sigma _{t}$ sweeps
out the \textquotedblleft annular strip\textquotedblright\ $\cup _{t\in
\lbrack 1-\epsilon ,1]}\Sigma _{t}$, on which there are exactly $n$ \textit{%
continuous} $n$th roots, i.e., continuous functions $g_{1},\cdots ,g_{n}$ on 
$\cup _{t\in \lbrack 1-\epsilon ,1]}\Sigma _{t}$ such that $\left[ g_{j}(w)%
\right] ^{n}=w$. As a result, for each $t\in \lbrack 1-\epsilon ,1]$, $%
\Gamma _{t}$ has exactly $n$ branches, with the $j$th branch parametrized by%
\[
\theta \mapsto g_{j}\left( 1+te^{i\theta }\right) 
\]%
(where the range of $\theta $ depends continuously on $t$ as well). It is
then clear that the length of $\Gamma _{t}$ varies continuously with $t$.
For $t$ ranging over $[1,1+\epsilon ]$, similar consideration leads to the
same conclusion concerning length but results in a different partition of $%
\Gamma _{1}$ into branches.

\medskip

\noindent \textbf{Remark}\quad Concerning the polynomial $p_{n}(z):=z^{n}-1$%
, we mention two connections, one a simple fact while the other possibly
(but improbably) fiction. First, $\left\vert p_{n}(z)\right\vert =t$ iff $%
|z|^{2n}-2\func{Re}(z^{n})+1=t^{2}$, and the latter is a real polynomial
equation of degree $2n$ in $(\func{Re}z,\func{Im}z)$; therefore the level
curves of $p_{n}$ are real algebraic curves of degree $2n$. Second, for $%
n\geq 2$, the level-$1$ curve of $p_{n}$ is conjectured around the
mid-twentieth century by \cite{E-P-H} to be no shorter than the level-$1$
curve of \textit{any} monic $n$th-degree polynomial. This conjecture remains
undecided except for the case $n=2$, which is proven true by \cite{E-H} near
the end of the last century. The level-$1$ curve of $z^{2}-1$ is a
Bernoulli's lemniscate.$\hfill \blacksquare $

\subsection{Length in Relation to Curvature and Derivative}

With the same notation as in \S 3.1, let $f$ be a \textit{nonconstant}
analytic function on an open connected subset $\Omega $ of $%
\mathbb{C}
$.

The real function $\left\vert f\right\vert $ is \textit{not} differentiable
at and only at simple zeros of $f$, as shown in \cite{B-D-N}. Where it is
differentiable,%
\[
\left\vert \nabla \left\vert f\right\vert \right\vert =\left\vert f^{\prime
}\right\vert =\left\vert \nabla \func{Re}f\right\vert 
\]%
by the Cauchy-Riemann equations. The critical points of $\left\vert
f\right\vert $, including points at which $\left\vert f\right\vert $ is not
differentiable, are isolated and Lemma \ref{Lemma KEY} is applicable to $%
\left\vert f\right\vert $.\footnote{%
If one wishes to avoid dealing with nondifferentiable functions, one may
work with $\left\vert f\right\vert ^{2}$ instead. A change of variable will
result in the same formula as obtained by working with $\left\vert
f\right\vert $.} Given the good behavior of length of level curves of $f$
(as seen in \S 3.1), the theory developed in \S 2.2 applies to $\left\vert
f\right\vert $ \textit{verbatim}, as long as the curvature $\kappa $ of
level curves is integrable over a disc around any critical point. While the
integrability of $\kappa $ is to be established in \S 5, we state brief
versions of the results in Lemma \ref{Lemma KEY} and Theorem \ref{Theorem
MAIN}.

\begin{theorem}
\label{Theorem Complex}Let $a,b$ be values attained by $\left\vert
f\right\vert $ such that $D:=\left\{ z:\left\vert f(z)\right\vert \in
\lbrack a,b]\right\} $ is compact. Then,

\begin{description}
\item[(a)] the area of $D$ equals $\int_{a}^{b}\left( \int_{\left\{
z:|f(z)|=t\right\} }\left( 1/\left\vert f^{\prime }(z)\right\vert \right)
\left\vert dz\right\vert \right) dt\,$;

\item[(b)] $\int_{a}^{b}L(t)\,dt=\iint\nolimits_{D}\left\vert f^{\prime
}\right\vert dA\,$;

\item[(c)] $L(b)=L(a)+\iint\nolimits_{D}\kappa \,dA\,$.
\end{description}
\end{theorem}

\noindent \textbf{Remark.}\textit{\quad }Suppose that a component of the
level-$c$ curve of $f$ bounds an open Jordan domain $D\subset \Omega $. In
this Remark, let $L(t)$ denote the length of the level-$t$ curve of $f|_{D}$.

This situation is markedly simpler than its counterpart for real functions
discussed for (\ref{EQ Jordan Boundary}) in \S 2.2. The extrema principle
for analytic functions implies that $\left\vert f\right\vert $ never again
attains the value $c$ in $D$, that $f$ vanishes somewhere in $D$, and that $%
\left\{ \left\vert f(z)\right\vert :z\in D\right\} =[0,c)$.

If, for each $t\in (0,c)$, the level-$t$ set of $f$ in $D$ is a simple
closed convex curve (in which case $f^{-1}(0)\cap D$ is necessarily a
singleton $\{z_{0}\}$ and $f^{\prime }$ does not vanish on $D\smallsetminus
\{z_{0}\}$), then%
\[
L^{\prime }(t)=\int_{\left\{ z\in D:|f(z)|=t\right\} }\frac{\kappa }{%
|f^{\prime }(z)|}\left\vert dz\right\vert >0 
\]%
for $t\in (0,c)$.\footnote{%
This positivity assertion is due to two reasons: our orientation stipulation
subject to which the formula for $L^{\prime }$ holds, and semi-definiteness
of the sign of $\kappa $ because of the convexity assumption on the curve $%
|f(z)|=t$.} Hence, $L$ is strictly increasing and we have the estimate%
\[
L(c)>\frac{1}{c}\iint\nolimits_{D}\left\vert f^{\prime }\right\vert dA\,%
\text{,} 
\]%
which is optimal as exemplified by the functions $z^{n}$ on $D(0;1)$ for $%
n\in 
\mathbb{N}
$.$\hfill \blacksquare $

\medskip

In closing this section, we say a word about curvature. Since each regular
level of $f$ is an analytic curve, its Schwarz function in principle can
yield its curvature; see \cite[p.\thinspace 45]{Davis}. In practice, with
machine-aided computation, the curvature formula (\ref{EQ Curvature})
applied to $\left\vert f\right\vert ^{2}$ suffices. However, we note a
formula for $\kappa $ purely in terms of $f$ and $f^{\prime }$ (which,
despite the best efforts of the author, has not been found in the
literature), and we leave its derivation to the reader.

\begin{proposition}
The level curves of an analytic function $f$ has curvature%
\[
\kappa =\frac{|f^{\prime }|}{|f|}-\frac{\left\langle \nabla |f|,\nabla
|f^{\prime }|\right\rangle }{|f^{\prime }|^{2}}\,\text{.} 
\]
\end{proposition}

\section{Applications and Extensions}

In \S 4.1, we give some geometric applications, among which is an analytic
characterization of circles. In \S 4.2, we deduce some additional identities
relating curvature and gradient and then extend our results to level
surfaces.

\subsection{Curve Evolution:\ Two Instances}

We show that some simple instances of curve evolution in the plane can be
framed in our elementary setting and addressed by our lines of inquiry.
Although most results that we deduce are well-known, it is the simplicity of
our method in obtaining them that we wish to illustrate. As a by-product, we
obtain a characterization of circles in terms of existence of certain
functions.

Consider a $C^{2}$ function $f$ on an open planar domain. Suppose that the
interval $[a,b]$ consists only of regular values attained by $f$ and that $%
f^{-1}([a,b])$ is connected and compact. For each $t\in \lbrack a,b]$, $%
f^{-1}(t)$ is then a simple closed curve. At $P\in f^{-1}([a,b])$, let $%
\mathbf{n}(P)$ be the outward unit normal at $P$ of the (simple closed)
curve $f^{-1}(f(P))$ and define%
\[
\sigma (P)=\left\langle \mathbf{n}(P),\frac{\nabla f(P)}{\left\vert \nabla
f(P)\right\vert }\right\rangle \text{ .} 
\]%
Being continuous and integer-valued, $\sigma $ is constant on (the
connected) $f^{-1}([a,b])$. By our stipulation in \S 2.1, the level curves
comprising $f^{-1}([a,b])$ all receive positive (resp. negative) orientation
iff $\sigma =+1$ (resp. $-1$). It is in this context that we will apply our
results to curve evolution.

Let $A(t)$ be the area enclosed by $f^{-1}(t)$. By Lemma \ref{Lemma KEY} and
Theorem \ref{Theorem MAIN},%
\[
L(b)=L(a)+\int_{a}^{b}\left( \int_{f^{-1}(t)}\frac{\kappa }{\left\vert
\nabla f\right\vert }ds\right) dt\,\text{;\ }A(b)=A(a)+\sigma
\int_{a}^{b}\left( \int_{f^{-1}(t)}\frac{1}{\left\vert \nabla f\right\vert }%
ds\right) dt\text{\thinspace .} 
\]

\noindent \textbf{Example 1}\quad Parallel closed curves and a
characterization of circles

Let $\alpha :[0,1]\rightarrow 
\mathbb{R}
^{2}$ be a simple closed $C^{2}$ curve. For $t\geq 0$, define $\gamma
_{t}:[0,1]\rightarrow 
\mathbb{R}
^{2}$ to be $\alpha +t\mathbf{n}$, where $\mathbf{n}$ is the outward unit
normal field along $\alpha $; i.e., the point $\gamma _{t}(u)$ is obtained
by traveling from $\alpha (u)$ along $\mathbf{n}(u)$ for a distance of $t$
units. As long as $t$ is sufficiently small, $\gamma _{t}$ defines a \textit{%
simple} closed curve; if $\alpha $ is convex, then $\gamma _{t}$ will be
simple and convex for any $t>0$. In any case, for any $t$ such that $\gamma
_{\tau }$ remains a simple closed curve for all $\tau \in \lbrack 0,t]$, we
seek the length $L(t)$ of $\gamma _{t}$ and the area $A(t)$ bounded by $%
\gamma _{t}$.

Define a function $f$ on the exterior of $\alpha $ by letting $f$ assume the
value $t$ on the curve $\gamma _{t}$; i.e., $f\left( \gamma _{t}(u)\right)
:=t$. Then, $\gamma _{t}$ is the level-$t$ curve of $f$, $\nabla f\left(
\gamma _{t}(u)\right) $ points in the outward normal direction $\mathbf{n}%
(u) $, and $\gamma _{t}$ is given the positive orientation according to the
stipulation in \S 2.1. Now, $\left\vert \nabla f\right\vert $, being the
derivative of $f$ with respect to distance along $\mathbf{n}$, clearly
equals $1$ (as $f$ \textit{is} the normal distance from $\alpha $). Hence,%
\begin{equation}
L(t)=L(0)+\int_{0}^{t}\left( \int_{\gamma _{\tau }}\kappa \,ds\right) d\tau
=L(0)+\int_{0}^{t}2\pi \,d\tau =L(0)+2\pi t\,\text{,}  \label{EQ L(t)}
\end{equation}%
and%
\begin{eqnarray}
A(t) &=&A(0)+\int_{0}^{t}\left( \int_{\gamma _{\tau }}1\,ds\right) d\tau
=A(0)+\int_{0}^{t}L(\tau )\,d\tau  \nonumber \\
&=&A(0)+L(0)t+\pi t^{2}\,\text{,}  \label{EQ A(t)}
\end{eqnarray}%
which are well known and stated for convex curves in \cite[p.\thinspace 47]%
{do Carmo}.

Note by the way an interesting consequence of (\ref{EQ L(t)}) and (\ref{EQ
A(t)}). If a simple closed convex curve were to undergo \textit{inward}
parallel evolution, then $A(t)$ and $L(t)$ vanish at the same time iff $%
A(0)=L(0)^{2}/4\pi $, in which case the initial curve is the optimal curve
for the isoperimetric inequality, i.e., a circle. For a more general notion
of parallel curves in relation to the isoperimetric inequality, see \cite[%
pp.\thinspace 79-85]{Guggenheimer}.

Formula (\ref{EQ L(t)}) shows that the evolving parallel curves lengthen at
the rate of $2\pi $ per unit normal distance $\Delta t$. We can draw the
same conclusion under a somewhat weaker hypothesis.

\begin{proposition}
\label{Prop. Area-Perimeter}Suppose that $f^{-1}([a,b])$ is compact and
that, for each $t\in \lbrack a,b]$, $f^{-1}(t)$\ is a simple closed curve on
which $\left\vert \nabla f\right\vert $\ equals a positive constant $c(t)$.
Then,

\begin{description}
\item[(a)] $\left\vert \nabla (L\circ f)\right\vert \equiv 2\pi $\ on $%
f^{-1}([a,b])$;

\item[(b)] the area of $f^{-1}([a,b])$\ equals $\left\vert
L(b)^{2}-L(a)^{2}\right\vert /4\pi $.
\end{description}
\end{proposition}

\noindent \textit{Proof. }Letting $\ell $ denote $L\circ f$ and $P_{0}\in
f^{-1}([a,b])$ with $f(P_{0})=t_{0}$, we calculate:%
\begin{eqnarray*}
\left\vert \nabla \ell (P_{0})\right\vert &=&\left\vert L^{\prime
}(t_{0})\right\vert \cdot \left\vert \nabla f(P_{0})\right\vert =\left\vert
\int_{f^{-1}(t_{0})}\frac{\kappa }{\left\vert \nabla f\right\vert }%
ds\right\vert \cdot \left\vert \nabla f(P_{0})\right\vert \\
&=&\frac{1}{\left\vert \nabla f(P_{0})\right\vert }\left\vert
\int_{f^{-1}(t_{0})}\kappa ds\right\vert \cdot \left\vert \nabla
f(P_{0})\right\vert =2\pi \text{ ,}
\end{eqnarray*}%
proving Part (a).

For Part (b), recall that $\sigma :=\left\langle \mathbf{n},\nabla
f/\left\vert \nabla f\right\vert \right\rangle $ is constant on $%
f^{-1}([a,b])$. Noting our orientation stipulation (subject to which Theorem %
\ref{Theorem MAIN}(b) holds), we have%
\[
L^{\prime }(t)=\int_{f^{-1}(t)}\frac{\kappa }{\left\vert \nabla f\right\vert 
}ds=\frac{1}{c(t)}\int_{f^{-1}(t)}\kappa \,ds=\sigma \frac{2\pi }{c(t)}\text{
.} 
\]%
Thus, $L^{\prime }$ has a definite sign on $[a,b]$, i.e., that of $\sigma $,
and so $L^{\prime }\circ f=\sigma \left\vert L^{\prime }\circ f\right\vert $%
. Hence%
\[
\left( L^{\prime }\circ f\right) \left\vert \nabla f\right\vert =\sigma
\left\vert L^{\prime }\circ f\right\vert \left\vert \nabla f\right\vert
=\sigma \left\vert \nabla \left( L\circ f\right) \right\vert =2\pi \sigma \,%
\text{.} 
\]%
Now consider%
\[
I:=\iint_{f^{-1}[a,b]}\left( L^{\prime }\circ f\right) \left\vert \nabla
f\right\vert dA\,\text{.} 
\]%
On one hand,%
\[
I=\iint_{f^{-1}[a,b]}2\pi \sigma dA=2\pi \sigma \iint_{f^{-1}[a,b]}1dA\,%
\text{;} 
\]%
on the other hand, by Theorem \ref{Theorem MAIN}(a) (with $h=L^{\prime }$),%
\[
I=\int_{a}^{b}L^{\prime }(t)L(t)dt=\frac{1}{2}\left[ L(b)^{2}-L(a)^{2}\right]
\,\text{.} 
\]%
Comparing the two expressions for $I$ proves the identity.\hfill $%
\blacksquare $

\smallskip

In light of this result, we raise and address the following question.

\smallskip

\noindent \textbf{Problem\quad }\textit{Given a simple closed }$C^{2}$%
\textit{\ curve }$C$\textit{\ bounding a compact Jordan domain }$R$\textit{\
and a point }$P$\textit{\ interior to }$R$\textit{, does there exist a }$%
C^{2}$\textit{\ Morse function }$f:R\rightarrow \lbrack 0,1]$\textit{\ such
that }$f^{-1}(0)=\{P\}$\textit{, }$f^{-1}(1)=C$\textit{, and, for each }$%
t\in (0,1]$\textit{, }$f^{-1}(t)$\textit{\ is a simple closed curve on which 
}$\left\vert \nabla f\right\vert $\textit{\ equals a positive constant?}

\smallskip

\noindent \textit{Solution}.\textit{\ }It is perhaps plausible that $R$
would have to be \textquotedblleft highly symmetric\textquotedblright\ (i.e.
a disc) about $P$. This is indeed the case, and can be deduced from
Proposition \ref{Prop. Area-Perimeter}(b). If there were such a function $f$%
, then the area of $f^{-1}([\epsilon ,1])$\ would equal $\left[
L^{2}(1)-L^{2}(\epsilon )\right] /4\pi $. Letting $\epsilon \rightarrow 0$,
it follows that the area $\mathcal{A}(R)$ of $R$ would equal $\mathcal{L}%
(C)^{2}/4\pi $, which, unless $C$ is a circle, violates the isoperimetric
inequality. In conclusion, only when $C$ is a circle can such a function
exist; in other words, existence of such a function characterizes the
circular domain that supports it.\hfill $\blacksquare $

\medskip

\noindent \textbf{Example 2}\quad The curve-shortening flow

Begin with a simple closed convex $C^{2}$ curve $\gamma
_{0}:[0,1]\rightarrow 
\mathbb{R}
^{2}$ with nonvanishing curvature. Move each point $P$ on $\gamma _{0}$
following the inward normal with initial speed equal to the\ curvature of $%
\gamma _{0}$ at $P$. Once motion begins (i.e., after an infinitesimal period
of time), the moved points comprise a new curve and will then be moved with
a new speed equal to the curvature of the new curve, so on and so forth. (As
points with differing curvature moves with different speed, it is
geometrically apparent that the evolution has the tendency to uniformize
curvature.) Obviously, we are attempting a verbal description of a
one-parameter family of evolving curves $\gamma _{t}:[0,1]\rightarrow 
\mathbb{R}
^{2}$ governed by a differential equation. The velocity of the point $\gamma
_{t}(u)$ is on one hand $\partial \gamma _{t}(u)/\partial t$ by definition
and on the other hand $\partial ^{2}\gamma _{t}(u)/\partial s^{2}$ by
prescription (where $s$ is arc length along $\gamma _{t}$); i.e.,%
\begin{equation}
\frac{\partial \gamma _{t}}{\partial t}=\frac{\partial ^{2}\gamma _{t}}{%
\partial s^{2}}\,\text{.}  \label{EQ PDE}
\end{equation}%
We assume the existence, uniqueness, and smoothness of solution (over a time
interval $[0,t]$) to this equation\footnote{%
This matter, belonging to the realm of partial differential equations, is
highly nontrivial and will take us too far afield.}, as well as the
simplicity and convexity of any future curve $\gamma _{t}$; see \cite%
{Chou-Zhu} for a thorough account or \cite[Appendix B]{Brakke} for a brief
account. The solution is known as the \textit{curve-shortening flow} for $%
\gamma _{0}$.

Define a function $f$ by letting $f$ assume the value $t$ on the curve $%
\gamma _{t}$. Then, $\gamma _{t}$ is the level-$t$ curve of $f$, along which 
$\nabla f$ points in the inward normal direction $-\mathbf{n}$. According to
our stipulation in \S 2.1, $\gamma _{t}$ is given the negative orientation
and its curvature $\kappa \leq 0$. Now, $\left\vert \nabla f\right\vert $,
being the derivative of time $t$ with respect to distance along $-\mathbf{n}$%
, clearly equals the reciprocal of the speed of evolution, i.e.,%
\[
\left\vert \nabla f\right\vert =1/\left\vert \kappa \right\vert =-1/\kappa \,%
\text{.} 
\]%
Therefore,%
\begin{eqnarray*}
L(t) &=&L(0)-\int_{0}^{t}\left( \int_{\gamma _{\tau }}\kappa ^{2}ds\right)
d\tau \\
\text{and\quad }A(t) &=&A(0)-\int_{0}^{t}\left( \int_{\gamma _{\tau }}\frac{1%
}{-1/\kappa }ds\right) d\tau =A(0)-2\pi t\text{\thinspace ,}
\end{eqnarray*}%
from which it follows that the evolution ceases by the time $A(0)/2\pi $ and
that $L$ is a strictly decreasing function because $\int_{\gamma _{\tau
}}\kappa ^{2}ds>0$.

\smallskip

\noindent \textbf{Remark\quad }When $\gamma _{0}$ is not convex, the PDE (%
\ref{EQ PDE}) still well-defines the curve-shortening flow of $\gamma _{0}$.
(In this case, points on $\gamma _{t}$ where the curve is not convex move
outward.) Then the difficulty in applying our method lies in the
construction of $f$, because $\gamma _{t}$ may intersect $\gamma _{t^{\prime
}}$ whereas the level curves of $f$ should not intersect. If we wish, we may
disentangle the intersecting curves by constructing the surface%
\[
S:=\left\{ \left( \gamma _{t}(u),t\right) :u\in \lbrack 0,1],t\in \lbrack
0,T]\right\} \subset 
\mathbb{R}
^{3} 
\]%
and then define $f$ on $S$ by $\left( \gamma _{t}(u),t\right) \mapsto t$.
The level-$t$ set of $f$ is then the curve $\gamma _{t}$ lifted to altitude $%
t$. Our results, once generalized to functions on surfaces, may conceivably
address this situation.\hfill $\blacksquare $

\subsection{Extensions}

Before we mention extensions of our results to level (hyper)surfaces, we
first note two additional relations between the curvature of level curves of
a two-variable function and its gradient field.

\begin{proposition}
Assume Hypothesis (\dag ). Then,

\begin{description}
\item[(a)] $\iint\nolimits_{f^{-1}([a,b])}\kappa
f_{x}dA=0=\iint\nolimits_{f^{-1}([a,b])}\kappa f_{y}dA\,$;

\item[(b)] if, in addition, $[a,b]$ is free of critical values and $%
f^{-1}([a,b])$ is connected, then $\iint\nolimits_{f^{-1}([a,b])}\kappa
\left\vert \nabla f\right\vert dA=\pm 2\pi (b-a)$.
\end{description}
\end{proposition}

\noindent \textit{Proof.} For Part (a), let $g$ in Lemma \ref{Lemma KEY} be
the vector-valued function $\kappa \nabla f$. Then,%
\[
\iint_{f^{-1}([a,b])}\kappa \nabla f\,dA=\int_{a}^{b}\left(
\int_{f^{-1}(t)}\kappa \frac{\nabla f}{\left\vert \nabla f\right\vert }%
ds\right) dt\,\text{.} 
\]%
It suffices to note that%
\[
\int_{f^{-1}(t)}\kappa \frac{\nabla f}{\left\vert \nabla f\right\vert }%
ds=-\int_{f^{-1}(t)}\kappa \mathbf{N\,}ds=-\int_{f^{-1}(t)}\frac{d\mathbf{T}%
}{ds}\mathbf{\,}ds=\mathbf{0} 
\]%
where the line integral is taken over all the components of $f^{-1}(t)$.

In Part (b), $f^{-1}([a,b])$ is diffeomorphic to $f^{-1}(a)\times \lbrack
a,b]$ and $f^{-1}(a)$ is then also connected. Hence, for $t\in \lbrack a,b]$%
, $f^{-1}(t)$, being diffeomorphic to $f^{-1}(a)$, is a simple closed curve
and $\int_{f^{-1}(t)}\kappa \,ds=\pm 2\pi $ (with the same sign in force for
all $t$, as shown in \S 4.1). Letting $g=\kappa \left\vert \nabla
f\right\vert $ in Lemma \ref{Lemma KEY} proves the assertion.\hfill $%
\blacksquare $

\medskip

To extend these and earlier results to level surfaces, we recall a few
differential-geometric preliminaries. Let $M$ be a connected, oriented,
compact $C^{2}$ surface in $%
\mathbb{R}
^{3}$. An \textit{orientation} on $M$ is a choice of a continuous unit
normal field $\mathbf{N}$ on $M$; the map%
\[
G:M\rightarrow S^{2};\;P\mapsto \mathbf{N}(P) 
\]%
is known as the \textit{Gauss map}. For $P\in M$ and a unit vector $\mathbf{v%
}\in T_{P}M:=\mathbf{N}(P)^{\perp }$, take the normal section $\Gamma _{%
\mathbf{v}}$ of $M$ whose velocity at $P$ is $\mathbf{v}$.\footnote{%
The curve $\Gamma _{\mathbf{v}}$ is the instersection of $M$ with the plane
determined by $\mathbf{N}$ and $\mathbf{v}$.} The (signed) curvature $\kappa
_{\mathbf{v}}$ of the curve $\Gamma _{\mathbf{v}}$ at $P$ is defined by%
\[
\left. \frac{d^{2}\Gamma _{\mathbf{v}}}{ds^{2}}\right\vert _{P}=\kappa _{%
\mathbf{v}}\mathbf{N}(P)\mathbf{\,}\text{.} 
\]%
The \textit{mean curvature} $H$ at $P$ is defined\footnote{%
Depending on the approach taken, there are various ways to define $H$. We
adopt the one that is the easiest to motivate.} to be, as the term suggests,
the mean of $\kappa _{\mathbf{v}}$ for $\mathbf{v}$ varying over the unit
circle in $T_{P}M$; in explicit terms,%
\[
H(P):=\frac{1}{2\pi }\int_{\left\{ \mathbf{v}\in T_{P}M:\left\vert \mathbf{v}%
\right\vert =1\right\} }\kappa _{\mathbf{v}}ds\,\text{.} 
\]%
(Euler's theorem on $\kappa _{\mathbf{v}}$ then implies that $H(P)=(\kappa
_{1}+\kappa _{2})/2$, where $\kappa _{1}$ and $\kappa _{2}$ are the extrema
of $\kappa _{\mathbf{v}}$, known as the \textit{principal curvatures} at $P$%
.) The sign of $H(P)$ obviously depends on the orientation of $M$. The 
\textit{Gaussian curvature} $K$ is simpler to define:%
\[
K(P):=\det \left( dG_{P}:T_{P}M\rightarrow T_{\mathbf{N}(P)}S^{2}\right) 
\text{\thinspace ,} 
\]%
which can be interpreted as the local (signed) area expansion factor at $P$
of the Gauss map. Clearly, $K(P)$ is independent of orientation on $M$. (As
a matter of fact, the principal curvatures $\kappa _{1}$ and $\kappa _{2}$
are the eigenvalues of $dG_{P}$, $K(P)=\kappa _{1}\kappa _{2}$, and $%
H(P)=\left( \limfunc{Tr}dG_{P}\right) /2$.)

Several facts are relevant here. Concerning $K$, we have%
\begin{equation}
\iint\nolimits_{M}K\mathbf{N}\,d\sigma _{M}=\mathbf{0}\text{\quad and\quad }%
\iint\nolimits_{M}K\,d\sigma _{M}=2\pi \chi (M)\text{\thinspace ,}
\label{EQ Gaussian Curvature}
\end{equation}%
where $d\sigma _{M}$ is the surface area form on $M$ and $\chi (M)$ is the
Euler characteristic of $M$. The second identity is the well-known
Gauss-Bonnet Theorem, while the first identity follows from a routine
calculation with differential forms.\footnote{%
For detail, define the vector-valued $2$-form $\omega $ on $S^{2}$ by
letting $\omega =\limfunc{Id}_{S^{2}}d\sigma _{S^{2}}$. Then $G^{\ast
}\omega =K\mathbf{N}\,d\sigma _{M}$, as can be verified pointwise. Hence,%
\[
\int_{M}K\mathbf{N}\,d\sigma _{M}=\int_{M}G^{\ast }\omega =\deg G\cdot
\int_{S^{2}}\omega \text{.} 
\]%
But $\int_{S^{2}}\omega =\int_{S^{2}}\limfunc{Id}\nolimits_{S^{2}}d\sigma
_{S^{2}}=\mathbf{0}$ due to cancellation of antipodal contributions.}
Concerning $H$, it is a fact (as shown in \cite[p.\thinspace 142]{do Carmo-2}%
) that%
\[
H=-\frac{1}{2}\func{div}\mathbf{N\,}\text{.} 
\]

Returning to our context, assume for simplicity that

\begin{quote}
$f$ is a $C^{2}$ function on a connected open set $\Omega \subset 
\mathbb{R}
^{3}$, $[a,b]$ is an interval of \textit{regular} values attained by $f$,
and $f^{-1}([a,b])$ is connected and compact.
\end{quote}

\noindent Then, for $t_{1},t_{2}\in \lbrack a,b]$, $f^{-1}(t_{1})$ and $%
f^{-1}(t_{2})$ are diffeomorphic compact surfaces. We orient $f^{-1}(t)$ by
letting $\mathbf{N}=-\nabla f/\left\vert \nabla f\right\vert $. Then, $H$
and $K$, both meaningful on $f^{-1}(t)$, become functions on $f^{-1}([a,b])$
and have explicit formulae in terms of the partial derivatives of $f$, as
shown in \cite[p.\thinspace 204]{Spivak}. With our choice of $\mathbf{N}$,%
\[
H=-\frac{1}{2}\func{div}\mathbf{N}=\frac{1}{2}\func{div}\frac{\nabla f}{%
\left\vert \nabla f\right\vert }\,\text{.} 
\]

Lemma \ref{Lemma KEY} in the present context takes the form%
\begin{equation}
\iiint\nolimits_{f^{-1}([a,b])}gdV=\int_{a}^{b}\left(
\iint\nolimits_{f^{-1}(t)}\frac{g}{\left\vert \nabla f\right\vert }d\sigma
\right) dt\text{\thinspace ,}  \label{EQ Level Surface}
\end{equation}%
which yields a formula for the volume of $f^{-1}([a,b])$ upon letting $%
g\equiv 1$.

Applying (\ref{EQ Level Surface}) with suitable choices of $g$, we have, as
consequences of (\ref{EQ Gaussian Curvature}),%
\[
\iiint\nolimits_{f^{-1}([a,b])}K\nabla f\,dV=\mathbf{0}\quad \text{and\quad }%
\iiint\nolimits_{f^{-1}([a,b])}K\left\vert \nabla f\right\vert \,dV=2\pi
(b-a)\chi (f^{-1}(a))\,\text{.} 
\]%
Writing $A(t)$ for the surface area of $f^{-1}(t)$, we have, for $t\in (a,b)$%
,%
\[
\frac{d}{dt}\iiint\nolimits_{f^{-1}([a,t])}\left\vert \nabla f\right\vert
dV=A(t)=A(a)+2\iiint\nolimits_{f^{-1}([a,t])}H\,dV\,\text{,} 
\]%
with the latter equality equivalent to%
\[
A^{\prime }(t)=\iint\nolimits_{f^{-1}(t)}\frac{2H}{\left\vert \nabla
f\right\vert }d\sigma \,\text{.} 
\]%
The proofs for these statements parallel those given for their counterparts
on level curves; except for the integral identity concerning $K\left\vert
\nabla f\right\vert $, we may as well allow $[a,b]$ to have critical values.
Similar to the treatment of curve evolution in \S 4.1, certain instances of
surface evolution, such as parallel expansion and mean curvature flow of a
convex surface, can be treated by these results.

We conclude with a cursory mention of level hypersurfaces of multivariable
functions. Given the background laid, it suffices to note that, for a
connected oriented one-codimensional compact submanifold $M$ of $%
\mathbb{R}
^{n+1}$ with Gauss map $G:M\ni P\mapsto \mathbf{N}(P)\in S^{n}$, we have the
mean curvature\footnote{%
There is a more general notion of mean curvature\ $\mu _{r}:=p_{r}(\kappa
_{1},\cdots ,\kappa _{n})/\binom{n}{r}$ where $p_{r}$ is the $r$th-degree
elementary symmetric polynomial and $\kappa _{j}$'s are the principal
curvatures, i.e., the eigenvalues of $dG$. Then, $\mu _{1}=H$ and $\mu
_{n}=K $. The integral of $\mu _{r}$ has geometric significance as well, but
discussion of them will take us too far afield.} and the Gaussian curvature%
\[
H:=\frac{1}{n}\limfunc{Tr}dG\quad \text{and}\quad K:=\det dG 
\]%
satisfying%
\[
H=-\frac{1}{n}\func{div}\mathbf{N}\quad \text{and}\quad \int_{M}K\,d\sigma
_{M}=\frac{1}{2}\chi (M)\sigma (S^{n})=\frac{\pi ^{\left( n+1\right) /2}}{%
\Gamma \left( \frac{n+1}{2}\right) }\chi (M)\,\text{.} 
\]

\section{Integrability of Curvature and a Generalized Morse Condition}

Finally, we prove that, if $f$ is a Morse function or the square of the
modulus of an analytic function, then the curvature $\kappa $ of its level
curves is integrable over a small neighborhood of any critical point. In
fact, the proof is valid under a more general hypothesis.

We first review the definition of improper multiple integral in a limited
context that meets our purpose; cf. \cite[pp.\thinspace 221-223]{Buck} or 
\cite[pp.\thinspace 257-260]{Courant}. Let $D$ be a closed disc centered at $%
P$ and let $g$ be continuous on $D\smallsetminus \{P\}$. Consider $%
I_{j}:=\iint\nolimits_{D_{j}}g\,dA$ where $\{D_{j}\}$ is an approximating
sequence of $D\smallsetminus \{P\}$, i.e., an expanding sequence of Jordan
measurable closed subsets of $D\smallsetminus \{P\}$ such that any interior
point of $D\smallsetminus \{P\}$ is eventually interior to some $D_{j}$. If $%
I_{j}$ converges to a limit that is independent of the choice of the
approximating sequence $\{D_{j}\}$, then $\lim I_{j}$ is taken to be the
value of the improper Riemann integral $\iint\nolimits_{D}g\,dA$.

It is shown in \cite[pp.\thinspace 221-223]{Buck} and \cite[pp.\thinspace
257-260]{Courant} that $\iint\nolimits_{D}g\,dA$ converges iff $%
\iint\nolimits_{D}\left\vert g\right\vert dA$ converges. If $|g(X)|$ is
bounded by $1/\left\vert XP\right\vert $ for $X\in D$, a calculation in
polar coordinates shows the convergence of $\iint\nolimits_{D}\left\vert
g\right\vert dA$ and hence of $\iint\nolimits_{D}g\,dA$.

\begin{proposition}
\label{Prop. Integrability}Suppose that, near the origin $O$, $%
f(x,y)=f(O)+p(x,y)+o(r^{n})$ for some homogeneous polynomial $p$ of degree $%
n\geq 2$. If $O$ (necessarily a critical point of $f$) is the only critical
point of $p$ (and hence an isolated critical point of $f$), then the
improper Riemann integral $\iint\nolimits_{D}\kappa \,dA$ converges on a
sufficiently small disc $D$ centered at $O$.
\end{proposition}

\noindent \textit{Proof.} Take $D$ to be so small that $O$ is the only
critical point of $f$ in it, in which case $\kappa $ is continuous on $%
D^{\prime }:=D\smallsetminus \{O\}$. (Recall (\ref{EQ Curvature}) in \S 2.1,
a formula for $\kappa $.)

Under our hypothesis,%
\[
f_{x}=p_{x}+r^{n-1}\epsilon _{1}\text{\quad and\quad }f_{y}=p_{y}+r^{n-1}%
\epsilon _{2}\,\text{,} 
\]%
where $\epsilon _{j}\rightarrow 0$ as $r\rightarrow 0$. Using polar
coordinates $(r,\theta )$ in $D^{\prime }$, we may write $p_{x}(\mathbf{r}%
)=r^{n-1}\alpha _{1}(\theta )$ and $p_{y}(\mathbf{r})=r^{n-1}\alpha
_{2}(\theta )$, where, e.g., $\alpha _{1}(\theta )=p_{x}(\cos \theta ,\sin
\theta )$. Hence, on $D^{\prime }$,%
\[
\left\vert \nabla f(\mathbf{r})\right\vert ^{2}=r^{2n-2}\left[ \left( \alpha
_{1}(\theta )+\epsilon _{1}\right) ^{2}+\left( \alpha _{2}(\theta )+\epsilon
_{2}\right) ^{2}\right] \,\text{.} 
\]%
As $O$ is assumed to be the only critical point of $p$, $\alpha _{1}(\theta
)^{2}+\alpha _{2}(\theta )^{2}\neq 0$ for $\theta \in \lbrack 0,2\pi ]$.
Letting $m=\min \left( \alpha _{1}(\theta )^{2}+\alpha _{2}(\theta
)^{2}\right) $, then $m>0$ and we have, for sufficiently small $r$, that $%
\left\vert \nabla f(\mathbf{r})\right\vert ^{2}>\frac{1}{2}mr^{2n-2}$ and
hence%
\[
\left\vert \nabla f(\mathbf{r})\right\vert ^{3}>Cr^{3n-3} 
\]%
for some constant $C>0$. By an entirely similar analysis, we can show that,
for sufficiently small $r$, there is some constant $M$ such that%
\[
\left\vert f_{xx}f_{y}^{2}-2f_{xy}f_{x}f_{y}+f_{yy}f_{x}^{2}\right\vert (%
\mathbf{r})<Mr^{3n-4}\text{.} 
\]

By (\ref{EQ Curvature}), we have%
\[
\left\vert \kappa (\mathbf{r})\right\vert =\frac{\left\vert
f_{xx}f_{y}^{2}-2f_{xy}f_{x}f_{y}+f_{yy}f_{x}^{2}\right\vert (\mathbf{r})}{%
\left\vert \nabla f(\mathbf{r})\right\vert ^{3}}<\frac{M}{C}\frac{1}{r} 
\]%
for sufficiently small $r$. Hence, $\kappa $ is integrable on $D$.\hfill $%
\blacksquare $

\medskip

The hypothesis of Proposition \ref{Prop. Integrability} motivates the
following definition.

\medskip

\noindent \textbf{Definition\quad }Let $n\geq 2$. A $C^{n}$ function $f$ on
a planar domain is said to satisfy \textit{the generalized Morse condition
of degree} $n$ if, at any critical point $(x_{0},y_{0})$, its $n$th-degree
Taylor polynomial equals%
\[
f(x_{0},y_{0})+p(x-x_{0},y-y_{0}) 
\]%
for a \textit{homogeneous} $n$th-degree polynomial $p$ whose \textit{only}
critical point is the origin.

\medskip

Lastly, we show that, for an analytic function $f$, the function $\left\vert
f\right\vert ^{2}$ satisfies the generalized Morse condition.

\begin{lemma}
\label{Lemma Analytic-Morse}For an analytic function $f$, the function $%
\left\vert f\right\vert ^{2}$ satisfies the generalized Morse condition.
\end{lemma}

\noindent \textit{Proof.} Let $Q=\left\vert f\right\vert ^{2}$. Suppose,
without loss of generality, that $0$ is a critical point of $Q$, in which
case either $f(0)=0$ or $f^{\prime }(0)=0$. For $z$ near $0$,%
\[
f(z)=f(0)+z^{n}g(z) 
\]%
for a unique integer $n\geq 1$ and a unique analytic function $g$ with $%
g(0)\neq 0$.

If $f(0)=0$,%
\[
Q(z)=\left\vert z^{n}g(z)\right\vert ^{2}=z^{2n}\left\vert g(z)\right\vert
^{2}=(x^{2}+y^{2})^{n}\left\vert g(z)\right\vert ^{2} 
\]%
and therefore $(x^{2}+y^{2})^{n}$ is the lowest-degree term in the Taylor
expansion of $Q(x,y)$ (as a real function). Obviously, $(x^{2}+y^{2})^{n}$
has only one critical point at the origin, and hence $Q$ satisfies the
generalized Morse condition of degree $2n$.

If $f^{\prime }(0)=0$, then $n\geq 2$. Let $a=f(0)$ and $b=g(0)$. Then%
\begin{eqnarray*}
Q(z) &=&\left\vert a+z^{n}g(z)\right\vert ^{2}=\left( a+z^{n}g(z)\right)
\left( \overline{a}+\overline{z}^{n}\overline{g(z)}\right) \\
&=&|a|^{2}+2\func{Re}(\overline{a}z^{n}g(z))+\left\vert z\right\vert
^{2n}\left\vert g(z)\right\vert ^{2} \\
&=&|a|^{2}+\func{Re}(2\overline{a}bz^{n})+o(\left\vert z\right\vert ^{n})\,%
\text{.}
\end{eqnarray*}%
Let $p(x,y)=\func{Re}(2\overline{a}b(x+iy)^{n})$, a homogeneous real
polynomial of degree $n\geq 2$. To see that $\nabla p(x,y)=\mathbf{0}$ iff $%
x=0=y$, it suffices to note that, for any analytic function $w$, $\left\vert
\nabla \func{Re}w\right\vert =\left\vert w^{\prime }\right\vert $ (by
Cauchy-Riemann equations). Hence, $Q$ satisfies the generalized Morse
condition of degree $n$.\hfill $\blacksquare $

\medskip

At long last, Proposition \ref{Prop. Integrability} and Lemma \ref{Lemma
Analytic-Morse} imply the following.

\begin{corollary}
The curvature of level curves of an analytic function is integrable around
each of its critical points.
\end{corollary}

\bigskip

\noindent \textit{\small Mathematics Department, Illinois State University,
Normal, Illinois}

\noindent \texttt{pding@ilstu.edu}

\end{document}